\documentclass[12pt]{amsart}
\usepackage{amssymb}

\newcommand{\tPa}[6]{\bibitem{#1} {#2}, \emph{#3}, {#4}, to appear (#5 pages).}

\newcommand{\Bgp}{{\Z^\N}}

\long\def\forget#1\forgotten{}
\newcommand{\issuenumber}{24}
\newcommand{\issuemonth}{March}
\newcommand{\issueyear}{2008}

\newtheorem{issue}{Issue}

\theoremstyle{definition}

\theoremstyle{remark}

\newcommand{\ed}{
\newpage

\section{Unsolved problems from earlier issues}

\begin{issue}
Is $\binom{\Omega}{\Gamma}=\binom{\Omega}{\Tau}$?
\end{issue}

\begin{issue}
Is $\ufin(\cO,\Omega)=\sfin(\Gamma,\Omega)$?
And if not, does $\ufin(\cO,\Gamma)$ imply
$\sfin(\Gamma,\Omega)$?
\end{issue}

\stepcounter{issue}

\begin{issue}
Does $\sone(\Omega,\Tau)$ imply $\ufin(\Gamma,\Gamma)$?
\end{issue}

\begin{issue}
Is $\fp=\fp^*$? (See the definition of $\fp^*$ in that issue.)
\end{issue}

\begin{issue}
Does there exist (in ZFC) an uncountable set satisfying $\sfin(\B,\B)$?
\end{issue}

\stepcounter{issue}

\begin{issue}
Does $X \nin \NON(\M)$ and $Y\nin\mathsf{D}$ imply that
$X\cup Y\nin \COF(\M)$?
\end{issue}

\begin{issue}[CH]
Is $\split(\Lambda,\Lambda)$ preserved under finite unions?
\end{issue}

\begin{issue}
Is $\cov(\M)=\fo$? (See the definition of $\fo$ in that issue.)
\end{issue}

\begin{issue}
Does $\sone(\Gamma,\Gamma)$ always contain an element of cardinality $\fb$?
\end{issue}

\begin{issue}
Could there be a Baire metric space $M$ of weight $\aleph_1$ and a partition
$\mathcal{U}$ of $M$ into $\aleph_1$ meager sets where for each ${\mathcal U}'\subset\mathcal U$,
$\bigcup {\mathcal U}'$ has the Baire property in $M$?
\end{issue}

\stepcounter{issue} 

\begin{issue}
Does there exist (in ZFC) a set of reals $X$ of cardinality $\fd$ such that all
finite powers of $X$ have Menger's property $\sfin(\cO,\cO)$?
\end{issue}

\begin{issue}
Can a Borel non-$\sigma$-compact group be generated by a Hurewicz subspace?
\end{issue}

\begin{issue}[MA]
Is there an uncountable $X\sbst\R$ satisfying $\sone(\BO,\BG)$?
\end{issue}

\begin{issue}[CH]
Is there a totally imperfect $X$ satisfying $\ufin(\cO,\Gamma)$
that can be mapped continuously onto $\Cantor$?
\end{issue}

\begin{issue}[CH]
Is there a Hurewicz $X$ such that $X^2$ is Menger but not Hurewicz?
\end{issue}

\begin{issue}
Does the Pytkeev property of $C_p(X)$ imply that $X$ has Menger's property?
\end{issue}

\begin{issue}
Does every hereditarily Hurewicz space satisfy $\sone(\BG,\BG)$?
\end{issue}

\begin{issue}[CH]
Is there a Rothberger-bounded $G\le\Bgp$ such that $G^2$ is not Menger-bounded?
\end{issue}

\begin{issue}
Let $\cW$ be the van der Waerden ideal.
Are $\cW$-ultrafilters closed under products?
\end{issue}

\begin{issue}
Is the $\delta$-property equivalent to the $\gamma$-property $\binom{\Omega}{\Gamma}$?
\end{issue}

\stepcounter{issue}

\general\end{document}}

\newcommand{\Cantor}{{\{0,1\}^\N}}

\newcommand{\fb}{\mathfrak{b}}

\newcommand{\fd}{\mathfrak{d}}
\newcommand{\fp}{\mathfrak{p}}

\newcommand{\NON}{{\mathsf   {NON}}}
\newcommand{\COF}{{\mathsf   {COF}}}

\newcommand{\M}{\mathcal{M}}

\newcommand{\cov}{\mathsf{cov}}

\newcommand{\R}{\mathbb{R}}

\newcommand{\fo}{\mathfrak{od}}

\renewcommand{\split}{\mathsf{Split}}
\newcommand{\bq}{\begin{quote}}
\newcommand{\eq}{\end{quote}}
\newcommand{\cO}{\mathcal{O}}
\newcommand{\B}{\mathcal{B}}
\newcommand{\BG}{\B_\Gamma}

\newcommand{\BO}{\B_\Omega}

\newcommand{\sone}{\mathsf{S}_1}    \newcommand{\sfin}{\mathsf{S}_\mathrm{fin}}

\newcommand{\ufin}{\mathsf{U}_\mathrm{fin}}

\newcommand{\Union}{\bigcup}
\newcommand{\nin}{\not\in}

\newcommand{\cW}{\mathcal{W}}

\newcommand{\N}{\mathbb{N}}
\newcommand{\Z}{\mathbb{Z}}

\newcommand{\sbst}{\subseteq}
\newcommand{\by}[2]{\par\hfill\emph{#1}, #2}
\newcommand{\nby}[1]{\par\hfill\emph{#1}}
\newcommand{\Tau}{\mathrm{T}}
\newcommand{\CE}{\textsc{CE}}

\newcommand{\be}{\begin{enumerate}}
\newcommand{\ee}{\end{enumerate}}
\newcommand{\bi}{\begin{itemize}}
\newcommand{\ei}{\end{itemize}}

\newcommand{\general}{\small\vfill\par\noindent\hrulefill\par
\noindent\textbf{Previous issues.} The previous issues of this
bulletin are available online at\\
\texttt{http://front.math.ucdavis.edu/search?\&t=\%22SPM+Bulletin\%22}
\\[0.1cm]
\textbf{Contributions.} Announcements, discussions, and open problems should be emailed
to \texttt{tsaban@math.biu.ac.il}\\[0.1cm]
\textbf{Subscription.}
To receive this bulletin (free) to your e-mailbox, e-mail us.
}

\newcommand{\arXiv}[5]{\subsection{#2}{#4}\par\hfill{\arx{#1}}\par\hfill\emph{#3}}

\newcommand{\AMS}[3]{\subsection{#1}~\par\hfill{\texttt{#3}}\par\hfill\emph{#2}}

\newcommand{\arx}[1]{\texttt{http://arxiv.org/abs/#1}}
\newcommand{\url}[1]{\bq\texttt{#1}\eq}
\newcommand{\online}[1]{The paper is available online at \url{#1}}

\title[$\mathcal{SPM}$ Bulletin \textbf{\issuenumber} (\issuemonth{} \issueyear)]{%
$\mathcal{SPM}$ Bulletin\\[0.5cm]
Issue number \issuenumber: \issuemonth{} \issueyear{} \CE{}}

\begin{document}
\maketitle

\tableofcontents

\section{Editor's note}

Enjoy.

\medskip

\by{Boaz Tsaban}{tsaban@math.biu.ac.il}

\hfill \texttt{http://www.cs.biu.ac.il/\~{}tsaban}

\section{Research announcements}

\subsection{A Wikipedia entry on topological games}
\url{http://en.wikipedia.org/wiki/Topological\_game}
\nby{Rastislav Telg\'arsky}

\AMS{On a fragment of the universal Baire property for $\Sigma^1_2$ sets}
{Stuart Zoble}
{http://www.ams.org/journal-getitem?pii=S0002-9939-08-08918-1}

\arXiv{0801.2132}
{The coarse classification of homogeneous ultra-metric spaces}
{Taras Banakh, Ihor Zarichnyy}
{We prove that two homogeneous ultra-metric spaces $X,Y$ are coarsely
equivalent if and only if $\mathrm{Ent}^\sharp(X)=\mathrm{Ent}^\sharp(Y)$ where
$\mathrm{Ent}^\sharp(X)$ is the so-called sharp entropy of $X$. This
classification implies that each homogeneous proper ultra-metric space is
coarsely equivalent to the anti-Cantor set $2^{<\omega}$. For the proof of
these results we develop a technique of towers which can have an independent
interest.}

\arXiv{0801.4723}
{Ramsey-like embeddings}
{Victoria Gitman}
{One of the numerous equivalent characterizations of a Ramsey cardinal
$\kappa$ involves the existence of certain types of elementary
embeddings for transitive sets of size $\kappa$ satisfying a large
fragment of ZFC. I introduce new large cardinal axioms
generalizing the Ramsey embeddings and show that they form a
natural hierarchy between weakly compact cardinals and measurable
cardinals.}

\arXiv{0801.4368}
{Proper and piecewise proper families of reals}
{Victoria Gitman}
{I introduced the notions of proper and piecewise proper families of reals to
make progress on an open question in the field of models of PA about whether
every Scott set is the standard system of a model of PA. A family of reals $X$ is
proper if it is arithmetically closed and the quotient Boolean algebra $X/fin$ is
a proper poset. A family is piecewise proper if it is the union of a chain of
proper families of size $\leq\omega_1$. Here, I explore the question of the
existence of proper and piecewise proper families of reals of different
cardinalities.}

\arXiv{0802.2705}
{Measures and their random reals}
{Jan Reimann, Theodore A. Slaman}
{We study the randomness properties of reals with respect to arbitrary
probability measures on Cantor space. We show that every non-recursive real is
non-trivially random with respect to some measure. The probability measures
constructed in the proof may have atoms. If one rules out the existence of
atoms, i.e.\ considers only continuous measures, it turns out that every
non-hyperarithmetical real is random for a continuous measure. On the other
hand, examples of reals not random for a continuous measure can be found
throughout the hyperarithmetical Turing degrees.}

\arXiv{0802.2550}
{Obtainable Sizes of Topologies on Finite Sets}
{Kari Ragnarsson and Bridget Eileen Tenner}
{We study the smallest possible number of points in a topological space having
$k$ open sets. Equivalently, this is the smallest possible number of elements in
a poset having $k$ order ideals. Using efficient algorithms for constructing a
topology with a prescribed size, we show that this number has a logarithmic
upper bound. We deduce that there exists a topology on $n$ points having $k$ open
sets, for all $k$ in an interval which is exponentially large in $n$. The
construction algorithms can be modified to produce topologies where the
smallest neighborhood of each point has a minimal size, and we give a range of
obtainable sizes for such topologies.}

\arXiv{0803.0676}
{Spaces of $\mathbb R$-places of rational function fields}
{Micha{\l} Machura and Katarzyna Osiak}
{In the paper an answer to the problem \emph{When different orders of R(X) (where R
is a real closed field) lead to the same real place ?} is given. We use this
result to show that the space of $\mathbb R$-places of the field
$\textbf{R}(Y)$ (where \textbf{R} is any real closure of $\mathbb R(X)$) is not
metrizable space. Thus the space $M(\mathbb R(X,Y))$ is not metrizable, too.}

\subsection{All properties in the Scheepers Diagram are linearly-$\sigma$-additive.}
In a present project, we study among other things the following preservation question.
Consider a property $P$ of topological spaces, which is not (provably) $\sigma$-additive.
Say that $P$ is \emph{linearly-$\sigma$-additive} if (provably), whenever
$X_1\sbst X_2\sbst\dots$, and each $X_n$ satisfies $P$, we have that
$\Union_n X_n$ satisfies $P$.

In the Scheepers Diagram, some of the properties
are $\sigma$-additive. The most interesting result for the properties which are not,
seems to be F. Jordan's recent theorem (see Issue 23) that $\sone(\Omega,\Gamma)$ is linearly-$\sigma$-additive.
Jordan's proof uses a combinatorial approach which is new to SPM. We supply an elementary proof, which seems
to provide the essence of Jordan's proof.
We also prove that $\sone(\Gamma,\Omega)$ and $\sfin(\Gamma,\Omega)$ are linearly-$\sigma$-additive.
This answers in the positive Problems 4.9 of \cite{SFT}, which asks whether these properties
are hereditary for $F_\sigma$ subsets.

Together with earlier results, we have that all properties in the Scheepers Diagram are linearly-$\sigma$-additive.

\nby{Tal Orenshtein and Boaz Tsaban}

\ed